\documentclass[10pt]{amsart}

\usepackage[active]{srcltx}
\usepackage{graphics,latexsym,dsfont,pictex,dcpic}
\usepackage{latexsym}
\usepackage{amsmath,amsthm,amssymb,bbm}
\usepackage{mathrsfs}

\newcommand{\ca}[1]{\mathcal{#1}}

\newcommand{\GF}{\mathds{F}}
\newcommand{\Naturals}{\mathds{N}}

\newcommand{\Rationals}{\mathds{Q}}

\newcommand{\Complex}{\mathds{C}}

\newcommand{\PSpc}{\mathbb{P}}

\theoremstyle{plain}

\theoremstyle{definition}

\numberwithin{equation}{thm}

\title[Disproof of Modularity]{
	Computational details on the disproof of modularity}

\author[Ralf Gerkmann]{Ralf Gerkmann}
\address{Universit\"{a}t Mainz, Fachbereich 17, Mathematik, 
	55099 Mainz, Germany}
\email{gerkmann@mathematik.uni-mainz.de}

\author[Mao Sheng]{Mao Sheng${}^\dagger$}
\address{East China Normal University, Dep.\ of Mathematics, 200062
	Shanghai, P.R.\ China}
\email{msheng@math.ecnu.edu.cn}
\thanks{${}^\dagger$ The second author was supported by the Chinese Science
	Foundation.}

\author[Kang Zuo]{Kang Zuo}
\address{Universit\"{a}t Mainz, Fachbereich 17, Mathematik, 55099 Mainz,
	Germany}
\email{kzuo@mathematik.uni-mainz.de}

\begin{document}

\setlength{\parindent}{0.0cm}
\fontsize{10}{14}\selectfont
\maketitle

\begin{abstract}
The purpose of these notes is to provide the details of the Jacobian ring
computations carried out in \cite{GMZ}, based on the computer algebra
system {\scshape Magma} \cite{Magma}.
\end{abstract}

\section{Finding an appropriate matrix}

Our first target is to construct the Jacobian ring associated to a Calabi-Yau 
$\tilde{X}$ as described in section 2 of \cite{GMZ}. This variety $\tilde{X}$
is defined by a matrix $A \in \Complex^{4 \times 8}$ subject to the condition
that all $4 \times 4$-minors are non-zero. In practice, we work over the field
$\Rationals$ or $\GF_p$, $p$ prime, instead of $\Complex$. The matrices
contained in the {\em hyperelliptic locus} were of special importance for us:
These are of the form $A := (a_{ij})$ with $a_{ij} = \lambda_j^i$, where 
$(\lambda_1,...,\lambda_8)$ denotes a tuple of $8$ distinct numbers in
$\Complex$.\medskip 

The {\scshape Magma} program below produces such an admissible matrix. There
are three options:

\begin{itemize}
\item[(i)] Work with a user-defined matrix. In this case one sets {\tt
randmat:=false} and {\tt hyperell:=false}. The line 
{\tt A:=RMatrixSpace(K,4,8)!0;} (marked with ``modify here'') has to be 
replaced by a statement that defines the desired matrix, e.g.\ 

\begin{center}
{\small
\begin{verbatim}

A:=Matrix(K,[
    [  1,  1,  1,  1,  1,  1,  1,  1],
    [  1,  2,  3,  4,  5,  6,  7,  8],
    [  1,  4,  9, 16, 25, 36, 49, 64],
    [  1,  8, 27, 64,125,216,343,512]]);

\end{verbatim}
}
\end{center}

In this case the non-degeneracy condition is {\em not} checked; it has to be
ensured by the user.

\item[(ii)] Use a random matrix. Set {\tt randmat:=true} and {\tt
hyperell:=false}. The program creates a random matrix in $K^{4 \times 8}$
which satisfies the non-degeneracy condition. The parameter {\tt randrange}
specifies the maximal absolute value of the matrix coefficients that are
created. 

\item[(iii)] Work with a user-defined matrix contained in the hyperelliptic
locus. In this case put {\tt hyperell:=true}, the value of {\tt randmat} is
irrelevant. In the line marked with ``modify here'' the definition  
{\tt la:=} can be replaced by any tuple of $8$ distinct elements in $K$.
\end{itemize}

For our implementation we used {\scshape Magma} version 2.11-2. It was carried
out under the operating system {\scshape Linux} on an ordinary personal
computer. The following program produces the desired matrix. 

\begin{center}
{\small
\begin{verbatim}

randmat     := true;
randrange   := 10;
hyperell    := true;

K := Rationals();
A := RMatrixSpace(K,4,8)!0; // <== modify here

genpos:=not(randmat) and not(hyperell);
while not(genpos) do
    N:=randrange;
    if hyperell then
        la:=[];
        for i:=1 to 8 do
            la[i]:=Random(2*N-2)-N+1;
        end for;
        la:=[1,2,3,4,5,6,7,8]; // <== modify here
        A:=[];
        for i:=0 to 3 do
            A[i+1]:=[];
            for j:=1 to 8 do
                A[i+1][j]:=la[j]^i;
            end for;
        end for;
    else
        A:=[];
        for i:=1 to 4 do
            A[i]:=[];
            for j:=1 to 8 do
                A[i][j]:=Random(2*N-2)-N+1;
            end for;
        end for;
    end if;
    A:=Matrix(K,A);
    seq:=[1,2,3,4,5,6,7,8];
    set:=SequenceToSet(seq);
    subsets:=Subsets(set,4);
    res:=true;
    for su in subsets do
        seq:=SetToSequence(su);
        B:=[];
        for i:=1 to 4 do
            B[i]:=[];
            for j:=1 to 4 do
                B[i][j]:=A[j][seq[i]];
            end for;
        end for;
        B:=Matrix(K,B);
        res:=res and not(Determinant(B) eq 0);
    end for;
    genpos:=res;
end while;

\end{verbatim}
}
\end{center}

In the first part of the program a random or hyperelliptic matrix
is created according to the user specification. The second part is to check
the non-degeneracy condition. This is done in the straightforward way: A loop
runs through all $4$-element subsets $S \subseteq \{ 1,...,8 \}$, creates the
submatrix that consists of the columns of $A$ given by the numbers in $S$, and
computes the determinant of this submatrix.

\section{Construction of the Jacobian ideal}

In order to carry out explicit computations in the cohomology of $\tilde{X}$ 
one replaces the intersection of four quadrics given by the matrix $A$ by a
certain toric hypersurface with defining equation
$$
F = y_1 f_1 + y_2 f_2 + y_3 f_3 + y_4 f_4 \qquad \text{with} \quad
f_1,...,f_4 \in K[x_0,...,x_7].
$$
as outlined in section 4 of \cite{GMZ}. The next step is to define the
polynomial ring and the Jacobian ideal associated to this hypersurface. 
The polynomial ring is created by

\begin{center}
{\small
\begin{verbatim}

PR:=PolynomialRing(K,12,"glex");
xn:=[];
for i:=0 to 7 do
    xn[i+1]:="x" cat IntegerToString(i);
end for;
yn:=[];
for i:=1 to 4 do
    yn[i]:="y" cat IntegerToString(i);
end for;
AssignNames(~PR,xn cat yn);
x:=function(i)
    return PR.(1+i);
end function;
y:=function(i)
    return PR.(8+i);
end function;

\end{verbatim}
}
\end{center}

The polynomial ring {\tt PR} has $12$ variables $x_0,...,x_7,y_1,...,y_4$ which
are accessed by {\tt x($i$)} and {\tt y($j$)} and diplayed as {\tt x$i$}, {\tt
y$j$}, respectively. Now the defining equation of the toric hypersurface can be
computed from the coefficient matrix $A$.

\begin{center}
{\small
\begin{verbatim}

f:=[];
for i:=1 to 4 do
    f[i]:=PR!0;
    for j:=0 to 7 do
        f[i]:=f[i]+A[i][j+1]*x(j)^2;
    end for;
end for;
F:=PR!0;
for i:=1 to 4 do
    F:=F+y(i)*f[i];
end for;

\end{verbatim}
}
\end{center}

The resulting equation is now stored in the variable {\tt F}. Finally, we
define the Jacobian ideal associated to this equation.

\begin{center}
{\small
\begin{verbatim}

I:=[];
for i:=1 to 12 do
    I:=Append(I,Derivative(F,i));
end for;
I:=Ideal(I);
I:=GroebnerBasis(I);
I:=IdealWithFixedBasis(I);

\end{verbatim}
}
\end{center}

\section{Deducing a cohomology basis}
\label{CohomologyBasis}

In this section we describe how to determine a basis for the cohomology of
$\tilde{X}$ as desribed in section 4 of \cite{GMZ}. In section 2 we defined an
action on the ring {\tt PR} by a certain finite group $N_1$. The cohomology of
$\tilde{X}$ is given by the $N_1$-fixed part of this ring modulo the Jacobian
ideal. It is a graded $K$-algebra 
$$
R = \bigoplus_{p=0}^3 R_p
$$
where the component $R_p$ is generated by the monomials with total degree
$2p$ in the variables $x_0,...,x_7$ and total degree $p$ in $y_1,...,y_4$. It
follows that a basis of $R$ can be computed by the following
procedure:

\begin{itemize}
\item[(i)] Enumerate all monomials of $x$-degree $2p$ and $y$-degree $p$.
\item[(ii)] Discard all monomials which do not agree with their normal form
with respect to the Jacobian ideal. This yields a basis of $R_p$.
\item[(iii)] Discard all monomials which are not fixed by the action of $N_1$.
The remaining monomials constitute a basis of $R_p$.
\end{itemize}

For (i) we define two additional polynomial rings with $8$ and $4$
indeterminates, respectively.

\begin{center}
{\small
\begin{verbatim}

XR:=PolynomialRing(K,8);
YR:=PolynomialRing(K,4);

\end{verbatim}
}
\end{center}

We also provide functions which map the indeterminates of {\tt XR} to
$x_0,...,x_7$ and the indeterminates {\tt YR} to $y_1,...,y_4$, respectively.
Using these functions, we are now able to enumerate all monomials of $x$-degree
$i$ and $y$-degree $j$ as follows: First, we use the {\scshape Magma} function 
{\tt MonomialsOfDegree} in order to produce all monomials of degree $i$ in {\tt
XR} and degree $j$ in {\tt YR}. The functions below cast these monomials into
elements of {\tt PR}, and the desired monomial set is obtained by multiplying
any element of the first set with any element of the second.

\begin{center}
{\small
\begin{verbatim}

xmon2mon:=function(xmon)
    res:=PR!1;
    for i:=0 to 7 do
        exp:=Degree(xmon,i+1);
        res:=res*(PR.(i+1))^exp;
    end for;
    return res;
end function;

ymon2mon:=function(ymon)
    res:=PR!1;
    for i:=1 to 4 do
        exp:=Degree(ymon,i);
        res:=res*(PR.(i+8))^exp;
    end for;
    return res;
end function;

\end{verbatim}
}
\end{center}

Step (ii) is based on the following simple function. 

\begin{center}
{\small
\begin{verbatim}

isbasiselt:=function(mon)
    nf:=NormalForm(mon,I);
    return (nf eq mon);
end function;

\end{verbatim}
}
\end{center}

Finally, step (iii) is carried out by the function below. It makes use of the
fact that the group $N_1$ is generated by $\sigma_a$, where $a$ runs through the
set
$$
\{ e_1+e_2,e_2+e_3,...,e_7+e_8,e_8+e_1 \}
$$ 
(see section 4 of \cite{GMZ} for the notation). 

\begin{center}
{\small
\begin{verbatim}

isHinvar:=function(mon)
    cond:=true;
    for i:=0 to 6 do
        deg1:=Degree(mon,i+1);
        deg2:=Degree(mon,i+2);
        cond:=cond and IsEven(deg1+deg2);
    end for;
    deg1:=Degree(mon,1);
    deg2:=Degree(mon,8);
    cond:=cond and IsEven(deg1+deg2);
    return cond;
end function;

\end{verbatim}
}
\end{center}

Based on these functions, the following loop now determines a basis for every
component the graded ring. The result is the nested array {\tt basis} in which 
{\tt basis[$p+1$]} contains the basis elements of $R_p$, given as elements
of the ring {\tt PR}.

\begin{center}
{\small
\begin{verbatim}

basis:=[];
for i:=0 to 2 do
    basis[i+1]:=[];
    xdeg:=2*i;ydeg:=i;
    xmons:=MonomialsOfDegree(XR,xdeg);
    ymons:=MonomialsOfDegree(YR,ydeg);
    for xmon in xmons do
        for ymon in ymons do
        mon:=xmon2mon(xmon)*ymon2mon(ymon);
            cond:=isHinvar(mon);
            cond:=cond and isbasiselt(mon);
            if cond then
                basis[i+1]:=Append(basis[i+1],mon);
            end if;
        end for;
    end for;
end for;
basis[4]:=[x(7)^6*y(4)^3];

\end{verbatim}
}
\end{center}

In order to save time, we define the component {\tt basis[$p+1$]} directly
instead of using the procedure. The reason is that the space of polynomials 
of bidegree $(6,3)$ is already quite large, so the computation would take
several seconds.

\section{The first characteristic subvariety}
\label{FirstCharSubv}

We briefly recall the definition of the characteristic subvariety $\ca{R}_1$
from section 3 in \cite{GMZ}. Let $e_1,...,e_9$ denote a basis of $R_1$ and
$e_1',...,e_9'$ a basis of $R_2$. Taking the multiplication map to its dual
yields a linear map
$$
\mu^* : R_2^* \longrightarrow S^2(R_1^*).
$$
Fix a bijection $\varphi : \{ (i,j) \in \Naturals^2 ~\mid~ 1 \leq i \leq j \leq
9 \} \longrightarrow \{ 1,...,45 \}$. Furthermore, let $A \in K^{9
\times 45}$ denote a representation matrix of $\mu^*$ with respect to
$(e_1')^*,...,(e_9')^*$ and $b_{\varphi(i,j)} := e_i^* e_j^*$. Then
$\ca{R}_1$ is the subvariety in $\PSpc^8$ given by the equations
$$
f_\ell := \sum_{i=1}^9 \sum_{j=i}^9 a_{\ell,\varphi(i,j)} z_i z_j \in
K[z_1,...,z_9].
$$
Thus in order to compute the variety, we have to carry out the following steps:

\begin{itemize}
\item[(i)] provide functions for the bijection $\varphi$ and its inverse
$\varphi^{-1}$
\item[(ii)] compute a representation matrix of the multiplication map $\mu :
S^2(R_1) \rightarrow R_2$ with respect to the basis
$\tilde{b}_{\varphi(i,j)} := e_i e_j$ and $e_1',...,e_9'$
\item[(iii)] obtain a representation matrix with respect to
$(e_1')^*,...,(e_9')^*$ and $\tilde{b}_{\varphi(i,j)}^* = (e_i e_j)^*$
by transposition
\item[(iv)] deduce matrix $A$ by base change from $(e_i e_j)^*$ to
$e_i^* e_j^*$
\item[(v)] write down the defining equations of $\ca{R}_1$
\end{itemize}

The bijections $\varphi$ and $\varphi^{-1}$ are provided by the following two
functions.

\begin{center}
{\small
\begin{verbatim}

r:=PolynomialRing(GF(2),9);
symm2order:=MonomialsOfDegree(r,2);
symm2order:=SetToSequence(symm2order);
Symm2Spc:=VectorSpace(K,#symm2order);

ijpos:=function(i,j)
    if i gt j then
        i_:=j; j_:=i;
    else
        i_:=i; j_:=j;
    end if;
    p:=9 - i_ + 1;
    offs:=45 -  1/2*p*(p+1);
    offs:=offs + (j_-i_) + 1;
    offs:=Integers()!offs;
    return offs;
end function;		

pos2ij:=function(pos)
    mon:=symm2order[pos];
    res:=[];
    for i:=1 to 9 do
        if Degree(mon,i) eq 2 then
            return [i,i];
        end if;
        if Degree(mon,i) eq 1 then
            res:=Append(res,i);
        end if;
    end for;
    return res;
end function;

\end{verbatim}
}
\end{center}

In order to carry out (ii) we need a function which computes, for an arbitrary
polynomial $g$ in the polynomial ring {\tt PR}, a coordinate vector of length
$1+9+9+1=20$ of the cohomology class represented by $g$ with respect to {\tt
basis}. 

\begin{center}
{\small
\begin{verbatim}

CSpc:=VectorSpace(K,20);
poly2vec:=function(poly)
    res:=CSpc!0;
    nf:=NormalForm(poly,I);
    offs:=0;
    for p:=0 to 3 do
        for i:=1 to #basis[p+1] do
            coe:=MonomialCoefficient(nf,basis[p+1][i]);
            res[offs+i]:=coe;
        end for;
        offs:=offs+#basis[p+1];
    end for;
    return res;
end function;

\end{verbatim}
}
\end{center}

The representation matrix is now obtained by computing the products of all
pairs of elements in {\tt basis[2]}. The $45$ rows of the matrix are given by
the $R_2$-parts of the coordinate vectors of these products.

\begin{center}
{\small
\begin{verbatim}

M:=[];
for l:=1 to 45 do
    M[l]:=[];
    ij:=pos2ij(l);i:=ij[1];j:=ij[2];
    mon1:=basis[2][i]; mon2:=basis[2][j];
    vec:=poly2vec(mon1*mon2);
    for k:=1 to 9 do
        M[l][k]:=vec[10+k];
    end for;
end for;
M:=Matrix(K,M);

\end{verbatim}
}
\end{center}

The implementation of the steps (iii) and (iv) is straightforward, 
provided one takes into account that $(e_i e_i)^* = e_i^* e_i^*$ 
and $(e_i e_j)^* = 2e_i^* e_j^*$. The final result is stored in the 
matrix named {\tt C}.

\begin{center}
{\small
\begin{verbatim}

N:=[];
for l:=1 to 45 do
    N[l]:=[];
    ij:=pos2ij(l);
    S:=SequenceToSet(ij);
    if #S eq 2 then
        la:=2;
    else
        la:=1;
    end if;
    for k:=1 to 9 do
        N[l][k]:=la*M[l][k];
    end for;
end for;
N:=Matrix(K,N);
C:=Transpose(N);

\end{verbatim}
}
\end{center}

Finally, the matrix {\tt C} is used to carry out (v). The resulting projective
scheme is stored in the variable {\tt charvar}.

\begin{center}
{\small
\begin{verbatim}

CSR:=PolynomialRing(K,9);
chareqs:=[];
for l:=1 to 9 do
    f:=CSR!0;
    for k:=1 to 45 do
        ij:=pos2ij(k);i:=ij[1];j:=ij[2];
        mon:=CSR.i*CSR.j;
        f:=f+C[l][k]*mon;
    end for;
    chareqs[l]:=f;
end for;
P8:=ProjectiveSpace(CSR);
charvar:=Scheme(P8,chareqs);

\end{verbatim}
}
\end{center}

The dimension and arithemtic genus are obtained by the function calls 
\begin{center}
{\tt Dimension(charvar)} \quad and \quad {\tt ArithmeticGenus(charvar)}.
\end{center}

\section{The second characteristic subvariety}

The computation of the characteristic subvariety $\ca{R}_2$ works essentially
in the same way as in the previous section. As before, let $e_1,...,e_9$ denote
a basis of $R_1$. By $e_1'$ we denote a vector which spans the $1$-dimensional
space $R_3$. The dual of the multiplication map yields
$$
\mu^* : R_3^* \longrightarrow S^3(R_1^*).
$$
In order to write down a representation matrix, we fix a bijection
$$
\ : \{ (i,j,k) \in \Naturals^3 ~\mid~ 1 \leq i \leq j \leq k \leq 
9 \} \longrightarrow \{ 1,...,210 \}.
$$
Let $A \in K^{1 \times 165}$ a representation matrix of $\mu^*$ with respect to
$e_1'$ and $b_{\varphi(i,j,k)} := e_i^* e_j^* e_k^*$. Then $\ca{R}_2$ is by
definition a hypersurface in $\PSpc^8$ given by the quation
$$
g := \sum_{i=1}^9 \sum_{j=i}^9 \sum_{k=j}^9 a_{1,\varphi(i,j,k)} z_i z_j z_k 
\in K[z_1,...,z_9].
$$
Thus we can proceed as follows.

\begin{itemize}
\item[(i)] provide functions for the bijection $\varphi$ and its inverse
$\varphi^{-1}$
\item[(ii)] compute a representation matrix of the multiplication map $\mu :
S^3(R_1) \rightarrow R_3$ with respect to the basis
$\tilde{b}_{\varphi(i,j,k)} := e_i e_j e_k$ and $e_1'$
\item[(iii)] obtain a representation matrix with respect to
$(e_1')^*$ and $\tilde{b}_{\varphi(i,j,k)}^* = (e_i e_j e_k)^*$
by transposition
\item[(iv)] deduce matrix $A$ by base change from $(e_i e_j e_k)^*$ to
$e_i^* e_j^* e_k^*$
\item[(v)] write down the defining equation of $\ca{R}_2$
\end{itemize}

Step (i) is achieved by

\begin{center}
{\small
\begin{verbatim}

r:=PolynomialRing(GF(2),9);
symm3order:=MonomialsOfDegree(r,3);
symm3order:=SetToSequence(symm3order);
Symm3Spc:=VectorSpace(K,#symm3order);

ijkpos:=function(i,j,k)
	S:=Sort([i,j,k]);
	i_:=S[1];j_:=S[2];k_:=S[3];
	offs:=0;
	for l:=1 to i_-1 do
		p:=10 - l;
		offs:=offs + 1/2*p*(p+1);
	end for;
	p0:=9-i_+1;
	j0:=j_-i_+1;
	k0:=k_-i_+1;
	offs2:=0;
	p0_:=p0 - j0 + 1;
	offs2:=1/2*p0*(p0+1) - 1/2*p0_*(p0_+1);
	offs2:=offs2 + (k0 - j0) + 1;
	offs:=offs+offs2;
	offs:=Integers()!offs;
	return offs;
end function;
	
pos2ijk:=function(pos)
	mon:=symm3order[pos];
	res:=[];
	for i:=1 to 9 do
		if Degree(mon,i) eq 3 then
			return [i,i,i];
		end if;
		if Degree(mon,i) eq 2 then
			res:=res cat [i,i];
		end if;
		if Degree(mon,i) eq 1 then
			res:=Append(res,i);
		end if;
	end for;
	return res;
end function;

\end{verbatim}
}
\end{center}

For step (ii) we can use the function {\tt poly2vec} from the previous section.
The representation matrix of the multiplication map $\mu$ is then computed by

\begin{center}
{\small
\begin{verbatim}

M:=[];
for l:=1 to 165 do
	M[l]:=[];
	ijk:=pos2ijk(l);i:=ijk[1];j:=ijk[2];k:=ijk[3];
	mon1:=basis[2][i]; mon2:=basis[2][j];
	mon3:=basis[2][k];
	vec:=poly2vec(mon1*mon2*mon3);
	M[l][1]:=vec[20];
end for;
M:=Matrix(K,M);

\end{verbatim}
}
\end{center}

For step (iii) and (iv), notice that $(e_i^3)^* = (e_i^*)^3$, $(e_i^2 e_j)^* =
3(e_i^*)^2 e_j^*$ and $(e_i e_j e_k)^* = 6e_i^* e_j^* e_k^*$ for $i,j,k$
pairwise distinct. 

\begin{center}
{\small
\begin{verbatim}

N:=[];
for l:=1 to 165 do
	N[l]:=[];
	ijk:=pos2ijk(l);
	S:=SequenceToSet(ijk);
	case #S:
		when 1:
			la:=1;
		when 2:
			la:=3;
		when 3:
			la:=6;
	end case;
	N[l][1]:=la*M[l][1];
end for;
N:=Matrix(K,N);
C:=Transpose(N);

\end{verbatim}
}
\end{center}

Now the characteristic variety is obtained by 

\begin{center}
{\small
\begin{verbatim}

CSR:=PolynomialRing(K,9);
chareq:=CSR!0;
for l:=1 to 165 do
    ijk:=pos2ijk(l);
    i:=ijk[1];j:=ijk[2];k:=ijk[3];
    mon:=CSR.i*CSR.j*CSR.k;
    chareq:=chareq+C[1][l]*mon;
end for;
P8:=ProjectiveSpace(CSR);
charvar:=Scheme(P8,chareq);

\end{verbatim}
}
\end{center}

\section{Computation of the Higgs field}

For the implementation of the {\em Plethysm method} (see section 3 in
\cite{GMZ}) it is neccessary to compute a representation matrices $A_w \in
K^{20 \times 20}$ of the maps
$$
\mu_w : R \longrightarrow R, \qquad x \mapsto wx \qquad \text{for} \quad w \in
R_1
$$
with respect to the basis computed in section \ref{CohomologyBasis}. The idea
is simple: Multiply every basis element with some fixed $w \in R_1$, use the
function {\tt poly2vec} from section \ref{FirstCharSubv} in order to transform
the products into row vectors and form a matrix out of these rows. This task is
carried out by the following code. The result is an array called {\em
thetamats} whose elements are the representation matrices corresponding to the
$9$ basis vectors of $R_1$.

\begin{center}
{\small
\begin{verbatim}

thetamats:=[];
for j:=1 to 9 do
    the:=basis[2][j];
    thetamat:=[];
    for p:=0 to 3 do
        for i:=1 to #basis[p+1] do
            poly:=the*basis[p+1][i];
            vec:=poly2vec(poly);
            vec:=ElementToSequence(vec);
            thetamat:=Append(thetamat,vec);
        end for;
    end for;
    thetamat:=Matrix(K,thetamat);
    thetamats:=Append(thetamats,thetamat);
    delete thetamat;
end for;

\end{verbatim}
}
\end{center}

\section{The induced Higgs field on symmetric $2$-space}

Each map $\mu_w : R \rightarrow R$ induces an endomorphism on $S^2(R)$. The aim
of this section is to represent the elements of $S^2(R)$ by vectors in
$K^{210}$, and to compute the action of $\mu_w$ with respect to this
representation. For the first part, we only have to fix a bijection
$$
\varphi : \{ (i,j) \in \Naturals^2 ~\mid~ 1 \leq i \leq j \leq 20 \}
\longrightarrow \{ 1,2,...,210 \}.
$$
If $b_1,...,b_{20}$ denotes a basis of $R$, then the isomorphism $S^2(R)
\rightarrow K^{210}$ is given by $b_i b_j \mapsto e_{\varphi(i,j)}$. The
following two functions provide such a bijection.

\begin{center}
{\small
\begin{verbatim}

r:=PolynomialRing(GF(2),20);
symm2order:=MonomialsOfDegree(r,2);
symm2order:=SetToSequence(symm2order);
Symm2Spc:=VectorSpace(K,#symm2order);

ijpos:=function(i,j)
    if i gt j then
        i_:=j; j_:=i;
    else
        i_:=i; j_:=j;
    end if;
    p:=20 - i_ + 1;
    offs:=210 -  1/2*p*(p+1);
    offs:=offs + (j_-i_) + 1;
    offs:=Integers()!offs;
    return offs;
end function;

pos2ij:=function(pos)
    mon:=symm2order[pos];
    res:=[];
    for i:=1 to 20 do
        if Degree(mon,i) eq 2 then
            return [i,i];
        end if;
        if Degree(mon,i) eq 1 then
            res:=Append(res,i);
        end if;
    end for;
    return res;
end function;

\end{verbatim}
}
\end{center}

Given an basis vector $w \in R_1$ and some $e_k \in K^{210} \cong S^2(R)$, our
next task is two compute the image $\mu_w(e_k)$. By the above isomorphism,
$e_k$ corresponds to an element of the form $b_i b_j$, and the induced map is
given by
$$
\mu(b_i b_j) = b_i \mu(b_j) + \mu(b_i) b_j.
$$
It means that we have to proceed as follows:

\begin{itemize}
\item[(i)] Use the matrices in {\tt thetamats} in order to compute the images
$\mu(b_i)$ and $\mu(b_j)$ as elements in $K^{20}$.
\item[(ii)] Compute the product $b_i \mu(b_j)$; if $\mu(b_j) =
\sum_{k=1}^{20} \alpha_k b_k$, then 
$$
b_i \mu(b_j) = \sum_{k=1}^{20} \alpha_k b_i b_k. 
$$
Use the bijection $\varphi$ in order to map this element to $K^{210}$. 
\item[(iii)] Do the same for $\mu(b_i) b_j$, and return the sum.
\end{itemize}

This is carried out by the following {\scshape Magma} code.

\begin{center}
{\small
\begin{verbatim}

symm2imthe:=function(th,pos)
    ij:=pos2ij(pos); i:=ij[1]; j:=ij[2];
    imi:=CSpc!0; imi[i]:=1;
    imj:=CSpc!0; imj[j]:=1;
    imj:=imj*thetamats[th];
    res1:=Symm2Spc!0;
    for k:=1 to 20 do
        if not(imj[k] eq 0) then
            res1[ijpos(i,k)]:=imi[i]*imj[k];
        end if;
    end for;
    imi:=CSpc!0; imi[i]:=1;
    imi:=imi*thetamats[th];
    imj:=CSpc!0; imj[j]:=1;
    res2:=Symm2Spc!0;
    for k:=1 to 20 do
        if not(imi[k] eq 0) then
            res2[ijpos(k,j)]:=imi[k]*imj[j];
        end if;
    end for;
    res:=res1+res2;
    return res;
end function;

\end{verbatim}
}
\end{center}

By linearity we can now compute the image of arbitrary elements in
$K^{210}$.

\begin{center}
{\small
\begin{verbatim}

symm2imthe_:=function(th,vec)
    res:=Symm2Spc!0;
    seq:=ElementToSequence(vec);
    for i:=1 to #seq do
        imvec:=symm2imthe(th,i);
        res:=res+vec[i]*imvec;
    end for;
    return res;
end function;

\end{verbatim}
}
\end{center}

For the application of the plethysm method, it is essential to carry out the
following task: Given a subspace $U \leq S^2(R)$, compute the sum
$$
\mu_{w_1}(U) \oplus \mu_{w_2}(U) \oplus \cdots \oplus \mu_{w_9}(U)
$$
of the images of all Higgs field maps. (Here $w_1,...,w_9$ denotes a basis of
the graded subspace $R_1$.) We use the built-in linear algebra functions of
{\scshape Magma} in order to provide this.

\begin{center}
{\small
\begin{verbatim}

symm2image:=function(U)
    B:=Basis(U);
    ims:=[];
    for bvec in B do
        for th:=1 to 9 do
            ims:=Append(ims,symm2imthe_(th,bvec));
        end for;
    end for;
    return sub<Symm2Spc|ims>;
end function;

\end{verbatim}
}
\end{center}

\section{The graded subspaces of symmetric $2$-space}

The grading $R = \oplus_{p=0}^3 R_p$ on the Jacobian ring induces a natural
grading 
$$
S^2(R) = \sum_{p=0}^6 S^2(R)_p
$$
on symmetric $2$-space. For the plethysm method, we have to compute the
iterated images of the graded subspaces $S^2(R)_p$ under the Higgs field
elements. In particular, for $0 \leq p \leq 6$ we have to determine the image
of $S^2(R)_p$ under the isomorphism $S^2(R)_p \cong K^{210}$. We provide a
function which returns for each $p$ a set $\{ i_1,...,i_{j_p} \}$ such that the
image of $S^2(R)_p$ in $K^{210}$ is spanned by the basis vectors
$e_{i_1},...,e_{i_{j_p}}$.\medskip

It seems to be the easiest strategy to compute the index set for each $p$ ``by
hand''. Let $b_1,...,b_{20}$ denote a basis of $R$. We assume that the indices
are chosen such that $R_0$ is spanned by $b_1$, $R_1$ is spanned
by $b_2,...,b_{10}$, $R_2$ by the elements $b_{11},...,b_{19}$ and finally $R_3$
by $b_{20}$. We first consider $p = 0$. In this case $S^2(R)_0$ is
generated by the image of $R_0 \otimes R_0$, which means that the image of
$S^2(R)_0$ in $K^{210}$ is spanned by $e_{\varphi(1,1)}$. The space $S^2(R)_1$
is generated by $R_0 \otimes R_1$, so the image is spanned by
$b_1 b_i = e_{\varphi(1,i)}$ with $2 \leq i \leq 9$. For $p=2$ we have to
consider the image of $(R_0 \otimes R_2) \oplus (R_1 \otimes R_1)$ etc. If we
work through all degrees up to $6$, we end up with the following
function.

\begin{center}
{\small
\begin{verbatim}

symm2degind:=function(p)
    case p:
        when 0:
            return [ijpos(1,1)];
        when 1:
            res:=[];
            for i:=1 to 9 do
                res:=Append(res,ijpos(1,1+i));
            end for;
            return res;
        when 2:
            res:=[];
            for i:=1 to 9 do
                res:=Append(res,ijpos(1,10+i));
            end for;
            for i:=1 to 9 do
                for j:=i to 9 do
                    res:=Append(res,ijpos(1+i,1+j));
                end for;
            end for;
            return res;
        when 3:
            res:=[ijpos(1,20)];
            for i:=1 to 9 do
                for j:=1 to 9 do
                    res:=Append(res,ijpos(1+i,10+j));
                end for;
            end for;
            return res;
        when 4:
            res:=[];
            for i:=1 to 9 do
                res:=Append(res,ijpos(1+i,20));
            end for;
            for i:=1 to 9 do
                for j:=i to 9 do
                    res:=Append(res,ijpos(10+i,10+j));
                end for;
            end for;
            return res;
        when 5:
            res:=[];
            for i:=1 to 9 do
                res:=Append(res,ijpos(10+i,20));
            end for;
            return res;		
        when 6:
            res:=[ijpos(20,20)];
            return res;
    end case;	
    return -1;
end function;

\end{verbatim}
}
\end{center}

The indices provided by {\tt symm2degind} now admit the computation of the
$p$-th graded subspace of $S^2(R)$.

\begin{center}
{\small
\begin{verbatim}

symm2degspc:=function(p)
    basis:=[];
    inds:=symm2degind(p);
    for i:=1 to #inds do
        vec:=Symm2Spc!0;
        vec[inds[i]]:=1;
        basis:=Append(basis,vec);
    end for;
    return sub<Symm2Spc|basis>;
end function;

\end{verbatim}
}
\end{center}

The above function now have to be invoked in the following way in order to make
the
plethysm method work.

\begin{center}
{\small
\begin{verbatim}

U51:=symm2degspc(1);
U42:=symm2image(U51);
U33:=symm2image(U42);
Dimension(U33);

\end{verbatim}
}
\end{center}

The result is $78$, whereas the modularity assumption would have implied a
dimension $\leq 65$.

\end{document}